\documentclass[12pt]{article}
\usepackage{amsmath, amsthm, amssymb, amscd}
\usepackage{graphicx}

\newtheorem{thm}{Theorem}
\newtheorem{cor}[thm]{Corollary}
\newtheorem{lemma}{Lemma}[section]
\newcommand{\Proof}{\noindent {\bf Proof: \quad} }
\newcommand{\Remark}{\noindent {\bf Remark: \quad} }

\newcommand{\mcg}{\ensuremath{MCG (\Sigma_g) }}
\newcommand{\freed}{\ensuremath{Free (\mathcal D_g) }}
\newcommand{\rhorg}{\ensuremath{\rho_{r,g} }}

\newcommand{\barimpsi}{\ensuremath{\overline {Im ( \psi) }}}
\newcommand{\expo}[1] {e^ { 2 \pi i #1}}

\newcommand{\ireven}[1]{\ensuremath{ I ( #1)}}

\newcommand{\mueven}{\mu}
\newcommand{\vgeven}{V_g}

\newcommand{\sumeven}[2]{ \sum_{\stackrel{k = #1}{ \text{even}}}^{#2} }
\newcommand{\sumevenj}[2]{ \sum_{\stackrel{j= #1}{ \text{even}}}^{#2} }

\title{SO(3) quantum invariants are dense}
\author{Helen Wong}
\date{}
\begin{document}

\maketitle
\author

\begin{abstract}
We show that when $r \geq 5$ is prime, the SO(3) Witten-Reshetikhin-Turaev quantum invariants for three-manifolds at the level $r$ form a dense set in the complex plane.  This confirms a conjecture of Larsen and Wang.
\end{abstract}

\section{Introduction}

Relatively recent work of Freedman, Larsen, and Wang in \cite{FLW} and \cite{LW} sought to better understand the $SO(3)$ quantum representation of the mapping class group of a surface.   
Though they were originally motivated by applications to quantum computation (see for example \cite{FKLW}), their results also have applications to studying topological three-manifolds.  

Historically, the first quantum theories originated with the Jones representation of the braid group and of the related Jones polynomial.    In \cite{Witten}, Witten conceived of combining the Jones polynomial with Chern-Simons quantum theory to produce a three-manifold invariant, which was subsequently realized by Reshetikhin and Turaev in \cite{RT}.  Further study of these constructions led to what is now referred to as examples of $SU(2)$ and $SO(3)$ Topological Quantum Field Theories.   Though it encompasses more, for us a TQFT  consists basically of a projective unitary representation of the mapping class group of a surface and a related three-manifold invariant.  As the $SU(2)$ theory may be derived from the $SO(3)$ theory, we will focus on the simpler $SO(3)$ version.  We provide more description in Section \ref{theory}.

Freedman, Larsen, and Wang in \cite{FLW} showed that the $SO(3)$ quantum representation at level $r=5$ has a dense image, and the result was generalized to include all prime and odd values of $r \geq 5$ by Larsen and Wang in \cite{LW}.
This paper builds on these density results and confirms a conjecture appearing in \cite{LW} regarding the related $SO(3)$ quantum invariant.   

For a choice of level $r$ and a root of unity $A$, let $I(M)$ denote the corresponding $SO(3)$ quantum invariant for a closed, connected, and orientable manifold $M$. \begin{thm} \label{main}
When $r$ is chosen to be a prime integer and $r \geq 5$ and when $A= i e^{2 \pi i /4r}$, the corresponding set of values of the $SO(3)$ quantum invariants for closed, orientable, connected three-manifolds $\{I(M)\}$ is dense in $ \mathbb C$.
\end{thm}
Note that since there are only a countable number of compact three-dimensional manifolds up to diffeomorphism, the $SO(3)$ quantum invariants can take on only a countable set of values.  As the quantum invariant $I(M)$ provides a lower bound on Heegaard genus of $M$ (\cite{T}), this result implies that for any genus $g$, there are infinitely many three-manifolds with Heegaard genus at least $g$.

The author would like to take this opportunity to thank Qihou Liu for his encouragement and enthusiasm, and also to thank her thesis advisor, Andrew Casson, for his insightful comments and continued support.  


\section{SO(3) quantum theory} \label{theory}

Fix $r \geq 3$ to be an odd integer, and choose $A$ to be one of the roots of unity $\pm e^{\pm 2 \pi i /4r}$ or $ \pm i e^{\pm 2 \pi i /4r}$.   We establish notation and briefly review the pertinent $SO(3)$ quantum theory associated to $r$ and $A$.  We refer the reader to either Turaev's book (\cite{T}) or Lickorish's survey (\cite{LSurvey}) for details on the basic constructions.   

Let $\Sigma_g$ be a closed, orientable surface of genus $g$.  Let $\mathcal D_g$ be the set of positive Dehn twists along simple closed curves in $\Sigma_g$.  The Dehn-Lickorish Theorem states that any self-homeomorphism of $\Sigma_g$ may be written as a finite product of Dehn twists. 
Thus $\mathcal D_g$ generates the mapping class group $\mcg$, and there exists a surjective map  $\nu: \freed \twoheadrightarrow \mcg$, where $\freed$ is the free group generated by $\mathcal D_g$.  

The $SO(3)$ quantum theory associated to $r$ and $A$ assigns to each Dehn twist on $\Sigma_g$ a unitary map of some finite dimensional vector space $V_g$ over $\mathbb C$.
For our purposes, we will not require the entire machinery involved in  constructing $V_g$.  We defer to the expositions of Turaev and Lickorish and only state the pertinent results here.
 
The vector space $V_g$ central to the $SO(3)$ theory is in fact a Hilbert space, with an associated inner product $\langle \cdot, \cdot \rangle : V_g \times V_g \rightarrow \mathbb C$. 
Further, as we mentioned, every simple closed curve $\gamma$ on $\Sigma_g$ induces a unitary map $\rhorg (\gamma)$ on $V_g$ with respect to this inner product.  
Let $U(V_g)$ denote the set of unitary maps on $V_g$, so that $\rhorg(\gamma) \in U(V_g)$.
Because each simple closed curve on $\Sigma_g$ corresponds to a Dehn twist along $\gamma$, and because products of unitary maps remain unitary, the $\rhorg(\gamma)$ extends to give a well-defined homomorphism $\rhorg: \freed \rightarrow U(V_g)$. 

However, Roberts in \cite{R} showed that upon composition with $\nu$, $\rhorg$ descends to only a {\sl projective} unitary representation of $\mcg$.  We have a  commutative diagram as follows:
\begin{equation} \label{CD}
 \begin{CD}
 \freed @>  \rhorg >>  U(V_g)\\
@V \nu VV @VV \pi V\\
\mcg @> \widetilde{\rhorg} >>  PU(V_g)
\end{CD}
\end{equation}
Here, the map $\pi$ is the standard projection map.  We will discuss $\pi$ further in Section \ref{pi}.

Roberts (\cite{Rirr}) also showed that that when $r$ is not only odd but  also prime, $\widetilde\rhorg$ is an irreducible representation.  From this, Freedman, Larsen and Wang (\cite{FLW}) and Larsen and Wang (\cite {LW}) wield deep results in representation theory to establish that $\widetilde\rhorg$ must have dense image.
\begin{thm} (\cite{FLW}, \cite{LW}) \label{FLW}
Let $r \geq 5$ be a prime integer, and let $A$ be one of $\pm e^{2 \pi i /4r}, \pm i e^{2 \pi i /4r} $.  When $g \geq 2$,  the map 
$\widetilde \rhorg : \mcg \longrightarrow PU(V_g)$ has dense image.
\end{thm}
For the purpose of this paper, it will be eaiser to consider a restatement of Theorem \ref{FLW}, which follows directly from the commutativity of (\ref{CD}).
\begin{cor} (\cite{FLW}, \cite{LW}) \label{FLWcor}
Let $r \geq 5$ be a prime integer, and let  $A$ be one of $\pm e^{2 \pi i /4r}, \pm i e^{2 \pi i /4r} $.  When $g \geq 2$,  the map 
$\pi \circ \rhorg : \freed \longrightarrow PU(V_g)$ has dense image.
\end{cor}

{\bf Remark}.  In the original papers \cite{FLW} and \cite{LW}, only the root of unity $A= \pm i e^{2 \pi i /4r}$ is considered.  However the same proof also applies to $A=\pm e^{2 \pi i /4r}$, as  the Dehn twist map has the same set of eigenvalues for either choice of $A$.\vspace{3pt}

Our main result, Theorem \ref{main}, seeks to convert Freedman, Larsen, and Wang's result into one about the $SO(3)$ quantum invariants for three-manifolds.  

In general, the $SO(3)$ quantum invariant can be defined for any odd value of $r$ and a choice of complex number $A$ with $A^4$ a primitive $r$th root of unity.  
But when $A$ is particularly one of $\pm e^{\pm 2 \pi i /4r}$ or $ \pm i e^{\pm 2 \pi i /4r}$ (so that the quantum representation exists), then, as we describe momentarily,
there is a simple relation between the  $SO(3)$ quantum invariant and the quantum representation.  

We begin by considering mapping cylinders of homeomorphisms on $\Sigma_g$. 
Let $H_g$ be a handlebody of genus $g$.
To every $f \in \freed$, we associate its mapping cylinder $M_f= H_g \cup_{\nu(f)} H_g$, obtained  by gluing two identical copies of $H_g$ by the homeomorphism $\nu(f) \in \mcg$ along their boundaries.  
It is important to note that because every three-manifold has a Heegaard splitting, any three-manifold may be represented this way.

The corresponding $SO(3)$ quantum invariant $I(M_f)$ may be viewed as deriving from the action of $\rhorg(f) \in U(V_g)$ on a single vector $v_\emptyset \in V_g$, the so-called normalized vacuum vector. In particular, 
\[ \ireven{M_f} = \alpha_f   \cdot \mu^{1-g}  \cdot \langle \rhorg(f) v_\emptyset, v_\emptyset \rangle.\]
$\alpha_f$ is a complex term dependent on $f$ and which captures information about  the framing of $M_f$.  Its inclusion is necessary to ensure that $I(M_f)$ is indeed a three-manifold invariant.  Though we provide no further details here, we do note that \emph{ $\alpha_f$ is always a $4r$th root of unity}.
The factor $\mu^{1-g}$ is a normalization constant associated to $v_\emptyset$ in $V_g$.  Further, $\mu$ is chosen so that $I(S^3) = \mu$ and $\mu^2 = \frac{ (A^2 - A^{-2})^2 }{ -r}$.

In the following proof of Theorem \ref{main}, it will be useful to regard the maps in $U(V_g)$ as matrices.  As a matter of convention, we will always choose the normalized vacuum vector $v_\emptyset$ as the first basis element.   Then $\langle \rhorg(f) v_\emptyset, v_\emptyset \rangle$ corresponds to the $(1,1)$-entry of a unitary matrix, which we write as $\rhorg(f)_{1,1}$ and so that \[ \ireven{M_f} = \alpha_f   \cdot \mu^{1-g}  \cdot \rhorg(f)_{1,1}.\]
We will capitalize on this relationship between $I(M_f)$ and $\rhorg(f)$ to prove Theorem \ref{main}.

We finally mention that the $SO(3)$ quantum invariant enjoys many nice properties.  For example, it behaves well under connected sum and change in orientation: 
\[ I(M \# M') = I(M) \cdot I(M') \quad \text{ and } \quad I(\overline M) = \overline{ I(M)} \]
for any three-manifolds $M$ and $M'$.  These properties will  play prominently in the proof.

\section{Projection map} \label{pi}

We will need a few facts from group theory, which we present in a general setting.

Let $U(n)$ denote the space of $n \times n$ unitary matrices, and let $PU(n)$ be the corresponding projective unitary matrices.  Recall that $PU(n)$ is defined as the space of equivalence classes of unitary matrices under the relation:  $U \sim V$ iff $U= e^{2 \pi i \theta} V$ for some $\theta \in \mathbb R$.  
This identification induces a projection $\pi: U(n) \rightarrow PU(n)$. 

In the following, we will consider $U(n)$ as a metric space with the Euclidean metric, where $|U-V| = \sqrt{ \sum_{i,j} | U_{i,j} - V_{i,j} |^2 }$.  Correspondingly, we have $ | \pi(U) - \pi(V)| = \min_{\theta, \nu} | e^{2 \pi i \theta} U - e^{2 \pi i \nu}V |$ in $PU(n)$.  This definition makes $\pi$ a continuous map.   For a subset $X$ of a metric space $Y$, let $\overline{X}$ denote its closure in $Y$.

\begin{lemma} \label{proj} 
Let $G$ be a group and $\psi : G \rightarrow U(n)$ be a group homomorphism.
Then the following statements are equivalent: 
\begin{enumerate}
\item 		$\overline{ Im( \pi \circ \psi)} = PU(n)$. 
\item  	Given any unitary $U \in U(n)$ and $\varepsilon >0$,
		there exists $g \in G$ and $\xi \in \mathbb R$ 
		so that $| e^{2 \pi i \xi}  \ \psi(g) - U| < \varepsilon$.
\item 		Given any unitary $U \in U(n)$,  
		there exists $\xi \in \mathbb R$ so that $ e^{2 \pi i \xi} U \in \barimpsi$.
\end{enumerate}
\end{lemma}

\Proof
$(1 \Rightarrow 2)$  
Suppose $U \in U(n)$ so that $\pi(U) \in PU(n) = \overline{ Im( \pi \circ \psi)}$.  
Thus there exists $ g \in G$ so that $| \pi \circ \psi (g) - \pi(U)| < \varepsilon$.
From the definitions, it follow that there exist $ e^{2 \pi i \xi_1}$ and $ e^{2 \pi i \xi_2}$ 
such that $ |e^{2 \pi i \xi_1} \ \psi(g) - e^{2 \pi i \xi_1} \ U| < \varepsilon$. So $ |e^{2 \pi i (\xi_1 - \xi_2)} \ \psi(g) - U| < \varepsilon$.

$(2 \Rightarrow 3)$
Let  $g_k $ and $\xi_k$ be so that $ |e^{2 \pi i \xi_k} \ \psi(g_k) - U| < 1/k$.  So $ |\psi(g_k) - e^{-2 \pi i \xi_k} \ U| < 1/k$.
By the compactness of the unit circle, there must be a subsequence $\{\xi_{k_j}\}$ so that $e^{ -2 \pi i \xi_{k_j}} \rightarrow e^{ -2 \pi i \xi_{\,}}$ for some $\xi \in \mathbb R$.  It follows that $ \psi(g_{k_j} )  \rightarrow e^{-2 \pi i \xi} \ U$.
 
$(3 \Rightarrow 1)$
If $e^{-2 \pi i \xi} \ U \in \barimpsi$, then $\pi(U) \in \pi (\barimpsi)$.  So $PU(n) \subseteq \pi ( \barimpsi)$.
Because $\pi$ is continuous,  $\pi ( \barimpsi) \subseteq \overline{ Im( \pi \circ \psi)}$.  
$\Box$
\\

\begin{lemma} \label{normal}
If  $\overline{ Im( \pi \circ \psi)} = PU(n)$, then $\barimpsi$ is a normal subgroup of $U(n)$.  
\end{lemma}

\Proof
Firstly, $\barimpsi$ is a subgroup of $U(n)$ since $\psi$ is a homomorphism.  
To show it is a normal subgroup, consider some  $S \in \barimpsi$ and $U \in U(N)$. 
The third statement in Lemma \ref{proj} implies that there exists $\xi \in \mathbb R$ so that $ e^{2 \pi i \xi} U \in \barimpsi$.
Then, $ U^{-1} S U = (e^{- 2 \pi i \xi} U^{-1}) S (e^{2 \pi i \xi} U  ) \in \barimpsi$. $\Box$

\section{Proof of Theorem \ref{main}} 

Assume that $ r\geq 5$ prime.  Let $A = i \expo{/4r} = i  \cos(  \pi/ 2r) -  \sin(\pi/2r)$ or $A = \expo{/4r} =  \cos(  \pi/ 2r)  + i \sin(\pi/2r)$,  so that Theorem \ref{FLW} and Corollary \ref{FLWcor}  apply.  
Let $g$ be any number at least 2.

Recall the definition $I(M_f)= \alpha_f   \cdot \mu^{-g}  \cdot \rhorg(f)_{1,1}$.  
For either choice of $A$, we have that
$\mu^2 = \frac{ (A^2 - A^{-2})^2 }{ -r} = \frac{ (2 i \ sin (\pi/r) )^2}{-r}$.  So $ \mu = \frac{2 \sin(\pi/r)}{ \sqrt{r}} $. Observe that $ 0 < \mu <1$ when $ r \geq 5$. 
Also recall that
$\alpha_f$ is a $4r$th root of unity.

We divide the proof into five lemmas.  
Let   $ \{ I(M) \} \subset \mathbb C $ denote the set of values of $SO(3)$ quantum invariants for closed, orientable, connected three-manifolds.   
We begin by showing that
 $ \{ I(M) \}  \cap \mathbb R$ is dense on the real line.

The first lemma is due to Larsen and Wang.  We reproduce the proof for the convenience of the reader.
\begin{lemma} (Larsen-Wang, \cite{LW}) \label{LWlemma}
Given a real $a \geq 0 $ and $\varepsilon > 0$, there exists a 3-manifold $N$ so that $\big{|} \,|I(N)| - a \big{|} < \varepsilon $.
\end{lemma}
\Proof
Pick $ g \geq 2 \text{ so that } 0 \leq a \cdot \mu^{g-1} <1$. Let $U \in U(\vgeven) $ with its $(1,1)$-entry being exactly $a \cdot \mu^{g-1}$, which we notate as $U_{(1,1)} = a \, \mu^{g-1}$.  
We apply Corollary \ref{FLWcor} in conjuction with Lemma \ref{proj} to show that
there exists 
$ f \in Free (\mathcal D)$ and $\xi$ with the property
$     | e^{2 \pi i \xi} \, \rhorg(f) - U | < \varepsilon  \, \mu^{g-1}. $
Restricting to the $(1,1)$ entry,
and $| e^{2 \pi i \xi} \, \rhorg(f)_{(1,1)} - a \, \mu^{g-1} | < \varepsilon  \, \mu^{g-1}$. 
This gives that
$  \big{|} \, | \mu^{1-g}  \, \rhorg (f)_{(1,1)}| - a \big{|} < \varepsilon. $

Finally observe that $ |I(M_f)|  = |\mu^{1-g} \cdot \alpha_f  \cdot \rhorg (f)_{(1,1)}|
             = | \mu^{1-g} \cdot \rhorg (f)_{(1,1)}| $ since $\alpha_f$ is a root of unity.
$\Box$\\


\begin{lemma}  \label{Rlemma}
Given a real $a \geq 0 $ and $\varepsilon > 0$, there exists $N$ 
so that $ I(N) \in \mathbb R$ and $| I(N) - a| < \varepsilon $.
\end{lemma}
\Proof
Since the square function is continuous and $a \geq 0$, 
there exists $M$ so that $\big{|} \, |I(M)|^2 - a \big{|} < \varepsilon $.
Properties of the invariant imply 
\[ I(M \# \overline M) = I(M) \cdot I(\overline{M}) = I(M) \cdot \overline{I(M)} = |I(M)|^2. \] 
Now let $N = M \# \overline M$. 
$\Box$\\

\noindent We next attempt to find a sequence of manifolds $M_{h_n}$ whose invariant $I(M_{h_n})$ converges to some complex number $b e^{2 \pi i \theta'}$ with $\theta' \notin \mathbb Q$.  This is the key step, allowing us to use Kronecker's Theorem to assert density.  
%
\begin{lemma}   \label{11lemma}
There exists a sequence $h_n \in \freed$ so that $\rhorg( h_n) _{1,1} \rightarrow  e^{2\pi  i  \theta}$ for some $\theta \notin \mathbb Q$.
\end{lemma}
\Proof 
Fix some $\nu \notin \mathbb Q$ and let
$ U = \left( \begin{array}{cc|l} 
                      e^{2\pi  i  \nu}& 0  & 0\\
                                     0& 1  & 0\\
\hline
                                     0& 0   &\text{Id}
\end{array} \right)$,
where Id denotes the identity matrix of size $\dim(\vgeven) -2$. 
There exists some $\xi \in \mathbb R$ with 
$e^{2\pi  i  \xi}U \in \overline{ Im( \rhorg)}$.
We are done if  $\nu + \xi \notin \mathbb Q$.

Else, $\nu + \xi \in \mathbb Q$ implies $\xi \notin \mathbb Q$.  
Let $ P = \left( \begin{array}{cc|l} 
                      0&  1 & 0\\
                      1&  0  & 0\\
\hline
                      0& 0   &\text{Id} 
\end{array} \right).$ 
Corollary \ref{FLWcor} in conjunction with Lemma \ref{normal} shows that
 $\overline{ Im( \rhorg)} $ is a normal subgroup of $U(\vgeven)$.
Thus the matrix $P^{-1}( e^{2\pi  i  \xi}U )P \in \overline{ Im( \rhorg)}$.  Its $(1,1)$ entry is $e^{2\pi  i  \xi}$.
$\Box$\\

\begin{lemma}  \label{notinqlemma}
There exists a sequence of 3-manifolds $M_{h_n}$ so that $I(M_{h_n}) \rightarrow b  \cdot e^{2\pi  i  \theta'}$ for some $b > 1$ and $\theta' \notin \mathbb Q$.
\end{lemma}
\Proof
The previous step provides us with a $\theta \notin \mathbb Q$ and sequence $h_n$ so that 
$ | \rhorg(h_n)_{1,1} - e^{2\pi  i  \theta} |  < (1/n) \,  \mu^{g-1}$.
We multiply through by $ \mu^{1-g} \cdot \alpha_{h_n}$ to get that $ | I(M_{h_n}) - {\mueven}^{-g}  \, \alpha_{h_n}  \, e^{2\pi  i  \theta} |  < 1/n. $

Because the $\alpha_{h_n}$  is a $4r$th root of unity, it can take on only a finite number of values.  
Thus we may pass to a subsequence where all $\alpha_{h_{n_j}}$ are equal, say to $\alpha$.  Thus, 
$ | I(M_{h_{n_j}}) - \mu^{1-g}  \, \alpha  \, e^{2\pi  i  \theta} |  < 1/n_j .$ 
Set $b= \mu^{1-g} $ and $ e^{2\pi  i  \theta'} = \alpha   \, e^{2\pi  i  \theta}$. 
$\Box$\\
\begin{lemma}  
Given $z= a  \, e^{2 \pi i \nu} \in \mathbb C$ with $a \geq 0$, and given $\varepsilon >0$, there exists a 3-manifold $M$ so that $| I(M) - z| < \varepsilon $.
\end{lemma}
\Proof
If $z=a$ then we are done by Step \ref{Rlemma}.    If not, take $a>0$.  We use an $\varepsilon$-thirds argument. 

Firstly, Step \ref{notinqlemma} gives $b > 1$, $\theta' \notin \mathbb Q$ and a sequence $\{ M_{h_n} \}$ 
such that  $I(M_{h_n}) \rightarrow b  \cdot e^{2\pi  i  \theta'}$.
Because $\theta' \notin \mathbb Q$, 
Kronecker's Theorem (Theorem $438$ in Chapter $23$ of \cite{HW}) says that $\{ e^{2 \pi i \theta' k} \}_{k \in \mathbb N}$ is dense in the unit circle.  We can thus find a $k \in \mathbb N$ with 
\[ | e^{2 \pi i \theta' k} -e^{2 \pi i \nu}| < \frac{\varepsilon }{3 \, a}. \]
Step \ref{Rlemma} gives a manifold $N$ with real  $I(N)>0 $ so that
\[ | I(N) - a  \cdot b^{-k} | < \frac{ \varepsilon}{ 3 \, b^k}, \]
By continuity, there also exists a sufficiently large $m \in \mathbb N$ 
so that $I(M_{h_m})$ has the property
\[ | I(M_{h_m})^k - (b  \cdot e^{2 \pi i \theta'})^k | < \frac{ \varepsilon}{ 3  \, I(N)}. \]

Combining the above three inequalities with the triangle inequality, we show that the manifold
$M= N \, \#_k  \, M_{h_m}$ with  $I(M)= I(N) I(M_{h_m})^k$ has the desired property.  
Namely,
\begin{align*}
| I(N) I(M_{h_m})^k  -  a  \cdot e^{2 \pi i \nu}| 
& \leq | I(N) I(M_{h_m})^k  -   I(N) b^k  e^{2 \pi i \theta' k} |  + \\
&       \quad \qquad  | I(N) b^k  e^{2 \pi i \theta' k} - a e^{2 \pi i \theta' k} | +
       | a e^{2 \pi i \theta' k} -  a  \cdot e^{2 \pi i \nu}| \\
&=  I(N) \cdot | I(M_{h_m})^k - (b  \cdot e^{2 \pi i \theta'})^k | +\\
&       \quad \qquad   b^k \cdot | I(N) - a  \cdot b^{-k} | + 
       a \cdot | e^{2 \pi i \theta' k} -e^{2 \pi i \nu}|  \\
& \leq \frac{\varepsilon}{3} + \frac{\varepsilon}{3} + \frac{\varepsilon}{3}. 
\end{align*}
$\Box$

\Remark  The proof of Theorem \ref{main} can be shortened by demonstrating a single $M$ whose invariant has irrational angle, i.e. $I(M) = a e^{2 \pi i \theta}$ with $\theta \notin \mathbb Q$.   This would replace the limiting arguments in Lemmas \ref{11lemma} and \ref{notinqlemma}.
However, the computation for $I(M)$ is generally quite complicated and difficult to write in closed form.    At least for some values of $A$, we do have an explicit, computable example which we sketch here.  

Let $M = T^2 \widetilde{\times}_f S^1$ be the mapping torus with fiber a torus and monodromy given by a single Dehn twist along a meridian curve.  Alternatively, this manifold may be thought of as that which results from integral surgery in $S^3$ along the Borromean link with framings $0$, $0$, and $1$.  By appealing to the definition of $I(M)$ and its properties, it can be seen that 
\[ I(M) = \mu \left( \sumeven{0}{r-2} A^{-k^2-2k}  \left( \frac{A^{2k+2} - A^{-2k-2}}{A^2 - A^{-2}} \right)^2  \right) \left( \sumevenj{0}{r-2}  A^{-j^2 - 2j} \right). \]
When $A=  e^{ 2 \pi i /4r}$ and $r \equiv 3 \mod 4$, elementary manipulations and a standard theorem on evaluating Gauss sums (see for instance \cite{HW}) show that 
\[ I(M) =  \frac{ \mu}{2} \left( \frac{ -A^2}{A^2 - A^{-2}} \right) \left( r \; - \;  i \sqrt{r} \right). \]
The $\frac{\mu}{2}$ term is a real number.  The next factor will have its $4r$th power on the real line (it is related to the $\alpha_f$ described in Section \ref{theory}).  The last factor will not have any of its integer powers being a real number.



\begin{thebibliography}{99}

\bibitem{Blanchet} Blanchet, C., {\sl Invariants on three-manifolds with spin structure}. Comment. Math. Helv. 1992, {\bf 67}, no. 3, 406-427.


\bibitem{BHMV} Blanchet, C., Habegger, N., Masbaum, G., and Vogel, P., {\sl Topological quantum field theories derived from the Kauffman bracket}. Topology 1995, {\bf 34}, no. 4, 883--927.

\bibitem{FLW} Freedman, M., Larsen, M., and Wang, Z., {\sl The two-eigenvalue problem and density of Jones representation of braid groups}. Commun. Math. Phys. 2002, {\bf 228}, 177-199.

\bibitem{FKLW} Freedman, M.,  Kitaev, A., Larsen, M.,  and Wang, Z., {\sl Topological quantum computation}. 
Bull. Amer. Math. Soc. 2003, {\bf40}, 31-38. 

\bibitem{HW} Hardy, G. and Wright, E. An Introduction to the Theory of Numbers. Fifth edition.
Clarendon Press, Oxford Univerity Press, 1979.


\bibitem{Lorig} Lickorish, W.B.R., {\sl Skeins and handlebodies}. Pacific J. Math 1993, {\bf 159}, 337-350.

\bibitem{L2r} Lickorish, W. B. R., {\sl The skein method for three-manifold invariants}. J. Knot Theory Ramifications 1993, {\bf 2}, no. 2, 171-194.

\bibitem{LSurvey} Lickorish, W. B. R.,  {\sl Quantum invariants of 3-manifolds}.  Handbook of geometric topology,  707--734, North-Holland, Amsterdam, 2002.

\bibitem{LW} Larsen, M. and Wang, Z., {\sl Density of the SO(3) TQFT representation of mapping class groups}.  Commun. Math. Phys. 2005, {\bf 260}, 641-658.

\bibitem{MRcentral}  Masbaum, G., and  Roberts, J., {\sl On central extensions of mapping class groups}. Math. Ann. 1995, {\bf 302}, no. 1, 131-150. 

\bibitem{R} Roberts, J.,  {\sl Skeins and mapping class groups}.
Math. Proc. Cambridge Philos. Soc, 1994, {\bf 115}, 53-77.

\bibitem{Rirr} Roberts, J., {\sl Irreducibility of some quantum representations of mapping class groups}. Knots in Hellas 1998, Vol. 3 (Delphi). J. Knot Theory Ramifications 2001 {\bf 10}, no. 5, 763--767.


\bibitem{RT} Reshetikhin, R. and Turaev, V., {Invariants of 3-manifolds via link polynomials and quantum groups}.  Invent. Math. 1991, {\bf 103}, 547-597.

\bibitem{T} Turaev, V. G.,   Quantum invariants of knots and 3-manifolds. de Gruyter Studies in Mathematics,  {\bf 18} . Walter de Gruyter \& Co., Berlin, 1994. 

\bibitem{Witten} Witten, E., {\sl Quantum field theory and the Jones polynomial}.  Commun. Math. Phys. 1989, {\bf 121}, 351-399.

\end{thebibliography}
\end{document}